\documentclass[12pt,reqno]{amsart}
\usepackage[usenames]{color}
\usepackage{amsmath, amsfonts, epsfig, amssymb, graphics}
\usepackage[colorlinks=true,citecolor=cyan]{hyperref}
\usepackage[mathscr]{euscript}

\def\N{{\mathbb N}}
\def\Q{\mathbb{Q}}

\def\S{{\mathbb S}}

\newcommand{\Harm}[1][W]{\mathscr{H}_{#1}}
\newcommand{\Diag}[1][W]{\mathscr{D}_{#1}}

\def\scalar#1#2{{\langle#1,#2 \rangle}} 

\newcommand{\newatop}[2]{\genfrac{}{}{0pt}{}{#1}{#2}}
\newcommand{\qbinom}[2]{\genfrac{[}{]}{0pt}{}{#1}{#2}}
\newcommand{\charac}{\raise 2pt\hbox{$\chi$}}

\def\auteur#1{{\sc #1}}
\def\titreref#1{{\em #1}}
\def\vol#1{{\bf #1}}
\def\mbf#1{{\bf #1}}

\newcommand{\pref}[1]{{\rm (\ref{#1})}}
\def\defn#1{\bleu{\bf #1}}

\def\bleu{\textcolor{blue}}
\def\rouge{\textcolor{red}}

\def\Q{\mathbb{Q}}

\def\S{{\mathbb S}}
\def\N{{\mathbb N}}

\newtheorem{proposition}{\bleu{Proposition}}
\newtheorem{conjecture}{\bleu{Conjecture}}
\newtheorem*{conj_autres}{Conjecture}

\begin{document}

\title[\rouge{Preliminary Version}]{{\bleu{New Formulas and Conjectures for the Nabla Operator}}}
\author[F.~Bergeron]{Fran\c{c}ois Bergeron}
\maketitle
\begin{abstract} The operator nabla, introduced by Garsia and the author, plays a crucial role in many aspect of the study of diagonal harmonics. 
Besides giving several new formulas involving this operator, we show how one is lead to representation theoretic explanations for conjectures about the effect of this operator on Schur functions.
\end{abstract}
 \parskip=0pt

{ \setcounter{tocdepth}{1}\parskip=0pt\footnotesize \tableofcontents}
\parskip=8pt

\section{Introduction}
Denoting by $\Delta_n(\mbf{x})$ the Vandermonde determinant in the variables $\mbf{x}=x_1,\ldots,x_n$, 
we consider the spaces
\begin{equation}\label{diag_troue}
  \Diag[{n;f}]:={\mathcal L}_{\partial,E}[f(\mbf{x})\Delta_n(\mbf{x})],
\end{equation}
where $f(\mbf{x})$ is a symmetric polynomial in the $\mbf{x}$-variables.
Here, for a polynomial $p(\mbf{x})$, we denote by ${\mathcal L}_{\partial,E}[p(\mbf{x}]$ the smallest vector space\footnote{The related spaces ${\mathcal L}_{\partial}[f(\mbf{x})\Delta_n(\mbf{x})]$ (with no application of the $E$-operators) are studied in \cite{ligne}.} (say over $\Q$) that contains $p(\mbf{x})$  and is closed under taking partial derivatives, as well as applications of the operators
    \begin{equation}\label{defoper}
            \bleu{E_{\mbf{u}\mbf{v}}^{(k)}:= \sum_{i=1}^n u_i\partial {v_i}^k},
      \end{equation}
 where $\mbf{u}$ and $\mbf{v}$ stand for any two of $r$ sets of $n$ variables:
    \begin{displaymath}
      \bleu{ \mbf{x}=x_1,\ldots,x_n},\qquad \bleu{\mbf{y}=y_1,\ldots,y_n},\qquad  \bleu{\mbf{z}=z_1,\ldots,z_n},\qquad {\rm etc}.\end{displaymath}
 The symmetric group $\S_n$ acts on   $\Diag[{n;f}]$ by \defn{diagonal} permutation of the various variable sets  $\mbf{u}$, sending $u_i$ to $u_{\sigma(i)}$, for $\sigma\in \S_n$. Although we restrict most of our discussion to the case of two sets of variables, Proposition~\ref{rest_general} holds in the more general context. It is easy to show that $\Diag[n]$ is a subspace of the space of diagonal harmonics. It is in fact the whole space for one (classical) and two (see~\cite{haimanvanishing}) sets of variables, and it seems to be so in general. For more on the case of $3$ or more sets of variables, see~\cite{several,three}.
      
The spaces  $\Diag[{n;f}]$ are \defn{graded} by  the \defn{(vector) degree}
          \begin{equation}\label{defndeg}
              \bleu{ \deg(f):= (\deg_{\mbf{x}}(f),\deg_{\mbf{y}}(f),\deg_{\mbf{z}}(f), \ldots )},
            \end{equation}
 and the action of $\S_n$ respects this grading, thus $\Diag[{n;f}]$ is a graded $\S_n$-module.
Recall that, for a (vector-degree) graded $\S_n$-module 
    \begin{displaymath}
      \bleu{\mathscr{V}=\bigoplus_{\mbf{d}\in \N^r} \mathscr{V}_\mbf{d}},
  \end{displaymath}
the coefficients of the \defn{Frobenius characteristic} (here denoted by  $\bleu{\mathscr{V}(\mbf{w};\mbf{q})}$), in the Schur basis of the $S_\mu(\mbf{w})$, correspond to (graded) multiplicities of irreducibles in the $\S_n$-module $\mathscr{V}$. This is to say that 
    \begin{displaymath}
      \bleu{\mathscr{V}(\mbf{w};\mbf{q})=\sum_{\mbf{d}=(d_1,d_2,\ldots,d_r) } q_1^{d_1}q_2^{d_2}\cdots  q_r^{d_r} \mathscr{V}_\mbf{d}(\mbf{w})},
  \end{displaymath}
with
    \begin{displaymath}
      \bleu{\mathscr{V}_\mbf{d}(\mbf{w})=\sum_{\mu\vdash n} a_{\mu;\mbf{d}} S_\mu(\mbf{w}),}
      \end{displaymath} 
where $a_{\mu;\mbf{q}}$ is the multiplicity of the irreducible accounted for $S_\mu(\mbf{w})$ in the homogeneous component $\mathscr{V}_{\mbf{d}}$ of degree $\mbf{d}$. Observe that the fact that  $\mathscr{V}$ is graded over $\N^r$ is encoded by the number of variables considered in $\mbf{q}$. 

When $f(\mbf{x})$ is a Schur function $s_\mu(\mbf{x})$, we write $\Diag[{n;\mu}]$ for $\Diag[{n;f}]$. To clarify any possible resulting ambiguity, let us agree that for $f(\mbf{x})=s_0(\mbf{x})=1$, we denote by $\Diag[n]$ the resulting space.
Our first objective here is to relate the Frobenius $\Diag[{n;k}](\mbf{w};\mbf{q})$ to conjectures of~\cite{nabla} concerning an operator $\nabla$ on symmetric functions. This operator is defined below in terms of Macdonald operators. On the way, we also deduce new identities concerning this operator $\nabla$ and some of its generalizations.

 
\section{Macdonald polynomials and the $\nabla$ operator}\label{intro}

\subsection*{Macdonald polynomials}
To go on with our story, we need to recall some basic facts about Macdonald polynomials and operators for which they are common eigenfunctions.
As usual (see \cite{macdonald}), we denote by $\lambda\preceq\mu$ the \defn{dominance order} on partitions. This is the (partial) order characterized by the fact that
   $$\lambda_1 + \lambda_2 + \cdots +
    \lambda_k\leq \mu_1 + \mu_2 + \cdots +\mu_k,\qquad \hbox{for all}\qquad k\geq 1,$$
setting $\mu_i=0$ or $\lambda_i=0$, whenever $i$ is larger than the number of parts of the underlying partition. Also as usual,  $\mu'$ stands for the \defn{conjugate} partition of $\mu$. 

Recall that the \defn{integral form Macdonald polynomials} $H_\mu(\mbf{w};q,t)$, with $\mu$ a partition of $n$ (we write $\mu\vdash n$), expand in the Schur function basis as
\begin{equation}\label{mac_kostka}
     \bleu{H_\mu(\mbf{w};q,t)=\sum_{\lambda\vdash n}
       K_{\lambda,\mu}(q,t)\, S_\lambda(\mbf{w})}.
   \end{equation}
The coefficients $K_{\lambda,\mu}(q,t)$ are known as the $q,t$-\defn{Kostka polynomials}. They have been shown, in~\cite{haimanvanishing}, to have positive integer coefficients.   
They form a linear basis of the ring $\Lambda$, of symmetric functions in the variables $\mbf{w}$, and   are characterized by the equations 
 \begin{equation}\label{eq_Hmu}
 \begin{array}{llll}
    \mathrm{(i)}\ \bleu{\displaystyle  \langle S_\lambda(\mbf{w}), H_\mu[(1-q)\,\mbf{w};q,t]\rangle=0},\qquad{\rm if}\qquad  \lambda\not\succeq\mu,\\[6pt]
    \mathrm{(ii)}\ \bleu{\displaystyle   \langle S_\lambda(\mbf{w}), H_\mu[(1-t)\,\mbf{w};q,t]\rangle=0},\qquad{\rm if}\qquad   \lambda\not\succeq\mu',\ \mathrm{and}\\ [6pt]
    \mathrm{(iii)}\ \bleu{\displaystyle   \langle S_n(\mbf{w}), H_\mu(\mbf{w};q,t)\rangle=1},
      \end{array}
 \end{equation} 
 involving the usual ``Hall'' scalar product on symmetric functions (for which the Schur functions are orthonormal).
Stated otherwise, (iii) says that $H_\mu(\mbf{w};q,t)$ is normalized so that $K_{n,\mu}(q,t)=1$, for all $\mu$. Since the one-part partition $(n)$ is largest in dominance order, the equations in \pref{eq_Hmu} imply that 
\begin{equation}\label{hall_little}
  \bleu{H_n(\mbf{w};q,t)=H_n(\mbf{w};q)=e_n\left[\frac{\mbf{w}}{1-q}\right]\, \prod_{k=1}^n (1-q^k)},
\end{equation}
using plethystic substitution notation (see~\cite{livre} and section~\ref{operator}). 
Observe that we may drop the parameter $t$, since it plays no role here.
This is a special case of a Hall-Littlewood polynomial (see \cite{macdonald}), so that $H_n(\mbf{w};q)$ has coefficients
\begin{equation}\label{cocharge_hook}
    \bleu{K_{\lambda,n}(q)=\sum_{\lambda(\tau)=\mu} q^{\mathrm{coch}(\tau)} },
  \end{equation}
  where the sum is over the set of standard tableaux $\tau$, of \defn{shape} $\lambda(\tau)=\mu$, with $\mathrm{coch}(\tau)$ standing for the \defn{cocharge} statistic (see \cite{macdonald}) of $\tau$. 
For example, 
\begin{eqnarray*}
H_1(\mbf{w};q)  &=&S_1(\mbf{w}),\\
H_2(\mbf{w};q)  &=&S_2(\mbf{w})+q\, S_{{11}}(\mbf{w}),\\
H_3(\mbf{w};q)  &=&S_3(\mbf{w})+ \left( {q}^{2}+q \right) S_{{21}}(\mbf{w})+{q}^{3}S_{{111}}(\mbf{w}),\\
H_4(\mbf{w};q)  &=&S_4(\mbf{w})+ \left( {q}^{3}+{q}^{2} +q\right) S_{{31}}(\mbf{w})+ \left( {q}^{4}+{q}^{2} \right) S_{{22}}(\mbf{w})\\
    &&\qquad+ \left( {q}^{5}+{q}^{4}+{q}^{3} \right) S_{{211}}(\mbf{w})+{q}^{6}S_{{1111}}(\mbf{w}).
\end{eqnarray*}
It has been shown in~\cite{haimanhilb} that  
      \begin{displaymath}
        \bleu{H_\mu(\mbf{w};t,q)=\Harm[\mu](\mbf{w};q,t)},
   \end{displaymath}
  where  $\Harm[\mu]:={\mathcal L}_\partial[\Delta_\mu]$ stands for the  \defn{Garsia-Haiman module}, obtained as the linear span of all partial derivatives of the  determinant
        \begin{displaymath}
        \bleu{\Delta_\mu:=\det\left(x_i^a y_i^b\right)_{\newatop{1\leq i\leq n}{(a,b)\in\mu}}},
   \end{displaymath}       
 with $(a,b)\in \mu$ meaning that $0\leq b\leq \ell(\mu)-1$, and $0\leq a\leq \mu_i-1$ (pairs being ordered lexicographically).
 
In our discussion, we  also make use of the $\S_n$-submodules (see~\cite{science_fiction}):
   \begin{displaymath}  
      \bleu{\Harm[\mu\downarrow]:=\bigcap_{\nu\rightarrow\mu} \Harm[\nu]},
   \end{displaymath}  
where $\nu\rightarrow\mu$ means that the partition $\nu$ precedes $\mu$ in the Young lattice.
 Using a set of heuristics discussed in \cite{science_fiction}, we can calculate $\Harm[\mu\downarrow](\mbf{w};q,t)$ in terms of Macdonald polynomials.

\subsection*{The $\nabla$ operator}\label{sec_nabla}
Since the $H_\mu(\mbf{w};q,t)$ form a basis of $\Lambda$, we can  define a linear operator $\nabla$  (read ``{nabla}'') on degree $n$ symmetric functions by  imposing that
     \begin{displaymath}  
         \bleu{\nabla(H_\mu(\mbf{w};q,t))=H_\mu(\mbf{w};q,t)\, \prod_{(a,b)\in \mu} q^at^b}.
      \end{displaymath}  
Thus the $H_\mu(\mbf{w};q,t)$ are eigenfunctions of $\nabla$ with simple, explicit eigenvalues expressed in term of the partition statistic 
    \begin{displaymath} 
         \bleu{ n(\mu):=\sum_i (i-1)\mu_i}.
    \end{displaymath}  
Indeed, the above definition is equivalent to ${\nabla(H_\mu(\mbf{w};q,t)):=q^{n(\mu')} t^{n(\mu)} H_\mu(\mbf{w};q,t)}$. 
 Recall that the main result in this context  is that
   \begin{equation}\label{formule_nabla} 
       \bleu{\Diag[n](\mbf{w};q,t)=\nabla(e_n(\mbf{w}))},
    \end{equation}
 where $\Diag[n]$ is being consider over two sets of $n$ variables, as expressed by the use of the two parameters $q$ and $t$ (rather than $q_1$ and $q_2$).

As discussed in~\cite{nabla}, explicit calculations reveal that the effect of $\nabla$ on Schur functions  is striking.
To expresse this, let us consider the matrix $\nabla^{(n)}:=(\nabla_{\lambda\mu})_{\lambda,\mu\vdash n}$, with
   \begin{displaymath}
       \bleu{\nabla_{\lambda\mu}(q,t):=\langle \nabla(S_\lambda(\mbf{w})),S_\mu(\mbf{w})\rangle}.
   \end{displaymath}
Here, we order partitions in decreasing lexicographic order, i.e.: $3$, $21$, $111$ for $n = 3$. The coefficients of these matrices are symmetric polynomials in the variables $q$ and $t$. We expand them as Schur polynomials (denoted with a lowercase ``s'' to further distinguish them from the $S_\mu(\mbf{w})$), omitting the variables $q$ and $t$.
 For $n=2$, we get 
$$ \nabla^{(2)} =\begin{pmatrix}0 & -s_{11} \\ 1 & s_1 \end{pmatrix}, $$
and for $n=3$,
$$  \nabla^{(3)} =\begin{pmatrix} 0 & s_{22} & s_{32} \\ 
0 & -s_{21} & -s_{31} \\ 
1 & s_2+s_1 & s_3+s_{11}\end{pmatrix}. $$ 
Inspection of these matrices, together with some theoretical considerations,  leads to 
\begin{conjecture}[see~\cite{nabla}]\label{nabla_conj}
For all $\lambda$ and $\mu$, 
 \begin{displaymath}
    \bleu{(-1)^{m(\lambda)} \nabla_{\lambda\mu}(q,t)}
 \end{displaymath}
has positive integer coefficients when expanded in terms of Schur polynomials,
with   
\begin{equation}
     \bleu{m(\lambda):=\binom{k}{2}+\sum_{\lambda'_i<(i-1)} (i-1-\lambda'_i)},
  \end{equation} 
 $k$ being the number of parts of $\lambda$.  
  \end{conjecture}
Similar conjectures have also been formulated in \cite{nabla} for iterates of $\nabla$, and for more general operators (see Section~\ref{operator}). Conjecture~\ref{nabla_conj} has been shown to hold in the special case $t=1$ by Lenart in \cite{lenart}. Among other results along these lines, let us also mention the work of Can and Loehr in \cite{can_loehr}.


\section{Representation theoretic conjecture for $\nabla(- S_{2,1^{n-2}})$}
Going back to the spaces introduced in section~\ref{intro}, we have the following.
\begin{proposition}\label{rest_general}
  The space $\Diag[{n;n}]$ is isomorphic to the restriction to $\S_n$ of the $\S_{n+1}$-module $\Diag[n+1]$.
 \end{proposition}
 \begin{proof}[\bf Proof.]
 Let us denote by $\mbf{x}'$ the set of $n+1$ variables $x_1,\ldots,x_n,x_{n+1}$, obtained by adding one extra variable to $\mbf{x}$.
We check that $\Diag[{n;n}]={\mathcal L}_{\partial,E}[e_n(\mbf{x})\Delta_n(\mbf{x})]$ is isomorphic to the restriction to $\S_n$ of $\Diag[n+1]={\mathcal L}_{\partial,E}[\Delta_{n+1}(\mbf{x}')]$ as follows. Considering its expansion as a polynomial in the variable $x_{n+1}$, it is easily seen that we may write  $\Delta_{n+1}(\mbf{x}')$ in the form
   \begin{equation}\label{expansion}
        \Delta_{n+1}(\mbf{x}')=e_n(\mbf{x})\Delta_n(\mbf{x}) + x_{n+1}\,g(\mbf{x}'),
   \end{equation}
 for some polynomial $g(\mbf{x}')$. Since $(\partial x_1+\ldots +\partial x_{n+1})\Delta_{n+1}(\mbf{x}')=0$, it is also easy to check that we may
 construct a basis of $\Diag[n+1]$ by application on $\Delta_{n+1}(\mbf{x}')$ of some operators built using only the $E$-operators and partial derivatives involving variables $x_k$, for which $1\leq k\leq n$. It follows that application of these same operators on the ``leading'' part $e_n(\mbf{x})\Delta_n(\mbf{x}) $ of \pref{expansion},  we get a basis of  $\Diag[{n;n}]$. Clearly this is all compatible with the action of $\S_n$ on the first $n$ variables in $\mbf{x}'$.
 \end{proof}
By general principles (see~\cite{livre}), we conclude that we have the  Frobenius equality
\begin{equation}\label{frob_restriction}
    \Diag[{n;n}](\mbf{w};q,t)=\partial p_1 \Diag[n+1](\mbf{w};q,t).
  \end{equation}
Using opeartors calculus outlined in section~\ref{operator}, it follows that 
\begin{proposition}\label{thm1}
In the case of two sets of variables, we have 
  \begin{equation}\label{reponse}
     \bleu{\Diag[{n;n}](\mbf{w};q,t)=\sum_{k=1}^n [k+1]_{q,t} \nabla( e_k(\mbf{w})\, e_{n-k}(\mbf{w}))},
  \end{equation}
  where
     \begin{eqnarray*}
           \bleu{[a]_{q,t}}&:=&\bleu{\frac{q^a-t^a}{q-t}}\\
                                   &=& \bleu{q^{a-1}+q^{a-2}t+\ldots + q t^{a-2}+t^{a-1}}.
     \end{eqnarray*}
  \end{proposition}
  Based on experiments and observations, we conjecture that more generally
\begin{conjecture}\label{conj111} For all $k$ between $1$ and $n-1$,
\begin{equation}\label{diag_n_k} 
       \bleu{\Diag[{n;k}](\mbf{w};q,t)=\sum_{j=0}^k
            [k-j+1]_{q,t}\,\nabla( e_j(\mbf{w})\,e_{n-j}(\mbf{w}))}.
     \end{equation}
 \end{conjecture}
Observe that, for $k=1$, (\ref{diag_n_k}) takes the form
     $$\Diag[{n ;1}](\mbf{w};q,t)=\nabla( e_j(\mbf{w})\,e_{n-j}(\mbf{w}))+(q+t)\,\nabla(e_n(\mbf{w})).$$
In view of the identity
 $S_{21^{n-1}}(\mbf{w})=e_1(\mbf{w})\,e_{n-1}(\mbf{w})-e_n(\mbf{w})$, we calculate directly that  Conjecture~\ref{conj111} implies that
\begin{equation}\label{frob_schur}
   \nabla(-S_{21^{n-1}}(\mbf{w}))= \Diag[{n ;1}](\mbf{w};q,t) -(q+t+1)\, \Diag[{n ;0}](\mbf{w};q,t).
 \end{equation}
Now, it is easy to check that
\begin{equation}\label{diag_inclusion}
    \Diag[n]\, \oplus\, e_1(\mbf{x})\,\Diag[n]\, \oplus\, e_1(\mbf{y})\,\Diag[n]\  \subseteq\ 
         \Diag[{n;1}].
   \end{equation}
This leads us to consider  the orthogonal complement $\mathcal{O}_Y $ (in the space $\Diag[{n;1}]$) of  the subspace corresponding to the left-hand-side of (\ref{diag_inclusion}).  The case $k=1$ of Conjecture~\ref{conj111} is thus seen to be equivalent to
\bleu{\begin{conjecture}
   $$   \mathcal{O}_Y (\mbf{w};q,t)=\nabla(-S_{21^{n-2}}(\mbf{w})).$$
  \end{conjecture}}
 Similar descriptions can be obtained for other values of $k$.
 A recent conjecture of N.~Loehr and G.~Warrington (see \cite{loehr_warrington})  describes in a combinatorial manner, for all partition $\lambda$, the expansion of $\nabla(S_\lambda(\mbf{w}))$ in terms of monomial symmetric functions. However, this gives only a purely enumerative description of the coefficients of the expansion of $\nabla(S_\lambda(\mbf{w}))$ in terms of the LLT symmetric polynomials (see \cite{LLT}).

\section{Operators}\label{operator}
We have already used plethystic substitution to describe $H_n(\mbf{w};q)$ in \pref{hall_little}. Recall that this operation turns symmetric functions into operators on the ring $\Lambda$, in the following manner. We first expand every symmetric function in terms of power sum $p_\lambda=p_{\lambda_1}\cdots p_{\lambda_k}$, and then apply the following rules:
 \begin{itemize}\itemsep=4pt
 \item $p_k[c]=c$ if $c$ is a constant;
 \item $p_k[x]=x^k$ if $x$ is a variable;
 \item $p_k[p_j]=p_{k\cdot j}$;
  \item $F[f+g]=F[f]+F[g]$, and $(F+G)[f]=F[f]+G[f]$;
  \item $F[f\cdot g]=F[f]\cdot F[g]$, and $(F\cdot G)[f]=F[f]\cdot G[f]$.
    \end{itemize}
  For the purpose of such calculations, $\mbf{w}$ is identified with $w_1+w+2+\ldots=p_1(\mbf{w})$.
As in \cite{nabla}, we can now introduce the operators $\mbf{D}_m$ (on symmetric functions $f(\mbf{w})$) defined by
the generating function identity 
\begin{equation}\label{def_D_m}
\sum_{m=-\infty}^{\infty} {\mathbf D}_m( f(\mbf{w}))\, \xi^m := f \!\left[\mbf{w} +{\alpha}/{\xi}\right]  \Omega'(\mbf{w};-\xi),
\end{equation}
with $\alpha=\alpha(q,t):=(1-q)(1-t)$, and
 $$\Omega'(\mbf{w};\xi):=\sum_{k\geq 0}e_k(\mbf{w})\,\xi^k=\prod_{i=1}^n {(1+r_i\,\xi)}.$$

Now, recall from \cite{nabla} that we have 
\begin{equation}\label{restriction_nabla}
   p_1^\perp\nabla=\alpha^{-1}\nabla\,{\mathbf D}_{-1}.
\end{equation}
For a symmetric function $f$, the operator $f^\perp$ is dual to multiplication by $f$, for the Hall scalar product.
We apply both sides of this operator identity to $e_{n+1}(\mbf{w})$, and calculate (directly using (\ref{def_D_m})) that 
\begin{eqnarray*}
      \partial p_1\nabla(e_{n+1}(\mbf{w}))&=&\nabla(\alpha^{-1}\,{\mathbf D}_{-1}(e_{n+1}(\mbf{w})))\\
             &=&\nabla\left(\sum_{k=0}^n [k+1]_{q,t}  e_k(\mbf{w})\, e_{n-k}(\mbf{w})\right)\\
              &=&\sum_{k=0}^n [k+1]_{q,t}  \nabla(e_k(\mbf{w})\, e_{n-k}(\mbf{w})).
  \end{eqnarray*}
Using \pref{formule_nabla}, we conclude that (\ref{reponse}) holds, since
\begin{eqnarray*}
    \Diag[{n;n}](\mbf{w};q,t)&=&\partial p_1\nabla(e_{n+1}(\mbf{w}))\\
                  &=&\sum_{k=0}^n [k+1]_{q,t}  \nabla(e_k(\mbf{w})\, e_{n-k}(\mbf{w})).
\end{eqnarray*}
This proves Proposition~\ref{thm1}.

Let us  recall some more operators identities of \cite{nabla} and \cite{explicit}.
Following \cite[(3.24)]{haglund} define the symmetric functions\footnote{Observe that we have had to change the notation used in \cite{haglund} from $E_{n,j}$  to $\varepsilon_{n,j}(\mbf{w};q)$. Otherwise, we would have had too many objects denoted by ``$E$''.} $\varepsilon_{n,j}(\mbf{w};q)$:
   $$\varepsilon_{n,j} (\mbf{w};q)=\sum_{k=0}^n F_k(\mbf{w};q)\, Z_k\Big|_{t_j},$$
 with 
   $$\sum_{k=0}^n F_k(\mbf{w};q)\, z^k=e_n\left[ \frac{\mbf{w}(1-z)}{1-q}\right],\qquad{\rm and}\qquad
      Z_k:=\sum_{i=0}^k(-1)^i q^{\binom{i+1}{2} -k\,i} (q;q)_i \qbinom{k}{i} t_i.$$
One can show (see \cite[Exercise 3.33]{haglund}) that
      \begin{equation}\label{epsilon1}
          (-q)^{n-1}\varepsilon_{n,1}(\mbf{w};q)=S_n(\mbf{w}).
     \end{equation}
It is also shown\footnote{Although without the use $\Harm[\mu\downarrow]$ which is tied with conjectures in \cite{science_fiction}.} in \cite[Formula (221)]{schroder} that
    \begin{equation}\label{epsilon_n_moins_1}
         q^{\binom{n}{2}}\, \varepsilon_{n,n-1}(\mbf{w};q)= -[n-1]_q \Harm[(n,1)\downarrow](\mbf{w}).
     \end{equation}
 One can further check (see \cite{haglund})  that 
  \begin{equation}\label{decomp_e}
     e_n(\mbf{w})=\sum_{j=1}^n \varepsilon_{n,j}(\mbf{w};q).
  \end{equation}
For any symmetric function $f$, as in \cite{nabla} we now consider the operators $\nabla_f$ having the
Macdonald polynomials $H_\mu(\mbf{w};q,t)$ as eigenfunctions with eigenvalue $f[B_\mu]$, writing $B_\mu$ for the polynomial $\sum_{(a,b)\in\mu} q^at^b$. Implicitly, a partition is here identified with the set of cells of its Ferrers diagram. 
In formula, the above definition can be stated as
\begin{equation}
   \nabla_f(H_\mu(\mbf{w};q,t)):=f[B_\mu]\, H_\mu(\mbf{w};q,t).
   \end{equation}
 Observe that we clearly have the operator identities $\nabla_{f+g}=\nabla_{f}+\nabla_{g}$ and $\nabla_{f\,g}=\nabla_{f}\nabla_{g}$. Moreover all these operators commute among themselves since they share a common basis of  eigenfunctions.
It has been shown in \cite{haimanvanishing} that, for all partitions $\mu$ and $\alpha$, one has the property
 \begin{equation}\label{nabla_schur_thm}
    \scalar{ \nabla_{S_\mu}\nabla(e_n(\mbf{w}))}{S_\alpha(\mbf{w})}  \in \N[q,t].
 \end{equation}
 In light of (\ref{decomp_e}) and computer calculations, the following refinement of Formula (\ref{nabla_schur_thm}) 
was also conjectured  to hold (see \cite[Conj.~4.12]{haglund}):
     \begin{equation}\label{haglund_epsilon}
              \scalar{ \nabla_{S_\mu}\nabla(\varepsilon_{n,k}(\mbf{w};q))}{S_\alpha(\mbf{w})}  \in \N[q,t], \qquad {\rm   for\ all\ }1\leq k\leq n.
     \end{equation} 
On the other hand, it was conjectured in \cite{nabla} that 
\bleu{\begin{conjecture} [BGHT 1999]\label{conj_BGHT}
 \begin{equation}
    \bleu{ \scalar{ \nabla_{S_\mu}(e_n(\mbf{w}))}{S_\alpha(\mbf{w})}  \in \N[q,t]},
 \end{equation}
 for all $\mu$ and $\alpha$.
\end{conjecture}}
For an analogous statement for the functions $\varepsilon_{n,k}(\mbf{w};q)$, see (\ref{nabla_epsilon}). 
 
We denote by $\nabla_k$ the operator obtained by choosing $f=e_k$ (the elementary symmetric function). In particular, when  restricted to the subspace of degree $n$ homogeneous symmetric functions, the operator $\nabla_{n}$ is simply the usual $\nabla$. In other terms, the usual $\nabla$ is 
    $$\nabla=\sum_k \nabla_k \pi_k,$$
 with $\pi_k$ denoting the projection of the ring of symmetric functions on its homogeneous component of degree $k$.
For the case $f=\sum_{k=0} e_k\,u^k$, we denote by $\Psi$ the operator $\Delta_f$. In other words, we have
\begin{equation}\label{def_psi}
     \Psi(H_\mu)= \prod_{(a,b)\in\mu} (1-q^at^b\,u)\,  H_\mu.
  \end{equation}
When restricted to degree $n$ symmetric functions, the operator $\Psi$ clearly expands as
\begin{equation}\label{etend}
    \Psi=\mathrm{Id}+u\,\nabla_1+u^2\,\nabla_2+\ldots +u^{n-1}\nabla_{n-1}+u^n\nabla_n .
 \end{equation}
For the sake of clarity, we now denote  by $\rho$   the operator $p_1^\perp$. 
Let us also write $\theta$ for the operator  $\alpha^{-1}{\mathbf D}_{-1}$, so that formulas appearing in formulas (I.12) of  \cite{nabla} may now be rewritten as
$$\begin{array}{clclc}
{\rm (a)} &  \mbf{D}_0= \mathrm{Id}-\alpha\,\nabla_1, &\qquad {\rm (b)} &  \alpha\, \mbf{D}_{k+1} =\mbf{D}_k\, \iota - \iota\, \mbf{D}_k ,\\[4pt]
{\rm (c)} &  \mbf{D}_1 = -\nabla \, \iota\, \nabla^{-1}, & \qquad {\rm (d)} &   \theta = \nabla^{-1} \, \rho\, \nabla,\\[4pt]
{\rm (e)} & \Psi^{-1}\,e_1 \,\Psi=e_1+u\,\theta, &\qquad  {\rm (f)} &  e_1\nabla_k = \nabla_k\,e_1+\nabla_{k-1}\theta.
\end{array}$$
We are now ready to derive more identities in the line of these.
\bleu{\begin{proposition}
We have the linear operators identities:
\begin{eqnarray}
     \Psi^{-1}\,\rho\,\Psi&=&\rho+u\,\theta,\label{maitre}\\
       \rho\nabla_k &=& \nabla_k\rho+\nabla_{k-1}\theta.\label{commute}
  \end{eqnarray}
 \end{proposition}}
 \begin{proof}[\bf Proof.]\ 
The first identity is easy to ``check'' by direct application on  the basis of Macdonald polynomials. 
Indeed, recalling the dual Pieri formula of Macdonald  
      \begin{displaymath}
            \rho\, H_\mu = \sum_{\nu\rightarrow \mu} c_{\mu,\nu}\, H_\nu,
        \end{displaymath}
  with $\nu\rightarrow\mu$ indicating that $\nu$ is covered by $\mu$ in the Young lattice of partitions, and the $c_{\mu,\nu}=c_{\mu,\nu}(q,t)$ being explicit rational fractions in $q$ and $t$ (see \cite{explicit}).
Applying the left hand side of (\ref{maitre}) to $H_\mu$ we get
      \begin{displaymath}
            \Psi^{-1}\,\rho\,\Psi (H_\mu)=\sum_{\nu\rightarrow \mu} c_{\mu,\nu} (1-q^{a_\nu} t^{b\nu}\,u) H_\nu,
        \end{displaymath}
where $(a_\nu,b_\nu)$ stands for the coordinates of the corner cell by which $\nu$ differs from $\mu$.  The right hand side of this last equality can clearly be written in the form
 $$ \Psi^{-1}\,\rho\,\Psi (H_\mu)=\rho\, H_\mu+u\,\nabla^{-1}\rho\nabla (H_\nu),$$
and we conclude using the operator equality $\theta=\nabla^{-1}\rho\nabla$ (see \cite[Formulas I.12]{nabla}). Comparing coefficients of $u^k$ in the expansion in (\ref{etend}), it follows from (\ref{maitre}) that (\ref{commute}) holds.
\end{proof}
 Iterative applications of (\ref{commute}) gives that
 \begin{equation}\label{mots}
    \bleu{\rho^n\nabla_{m} =\sum_{j=0}^{n}\sum_{\alpha\in {\mathcal A}_{j,n-j}} \nabla_{m-j} \,\alpha},
 \end{equation}   
with ${\mathcal A}_{j,k}$ denoting the set of $\binom{j+k}{k}$ operators that can be obtained by composing $j$ copies of $\rho$ and $k$ copies of $\theta$ in all possible order.
Applying both sides of (\ref{mots}) to $e_n(\mbf{w})$, we find that
 \begin{equation}\label{get_dim}
   \bleu{ \rho^n\nabla_{n-1}(e_n(\mbf{w})) = \sum_{k=0}^{n-1} \theta^k\rho\, \theta^{n-1-k} e_n(\mbf{w})},
 \end{equation}
 since all other terms of the resulting sum vanish, and $\nabla_0=\mathrm{Id}$. 
It happens (as we will see below) that $\nabla_{n-1}(e_n(\mbf{w}))$ corresponds to the bigraded Frobenius characteristic of some $\S_n$-module. Formula~(\ref{get_dim}) can be used to get an explicit formula for the Hilbert series of this module.
 
 Other exciting consequences of (\ref{maitre}) and (\ref{commute}) are that
    \bleu{ \begin{eqnarray}
   \rho \, \Psi(p_n(\mbf{w})) &=&(-1)^{n-1} u\,[n]_t[n]_q\,\Psi(e_n(\mbf{w})),\\
  \nabla_{n-1}(p_n(\mbf{w})) &=& \frac{[n]_t[n]_q}{q^{n-1}t^{n-1}}\nabla(h_n(\mbf{w})).
 \end{eqnarray}}
 It is also interesting that we have
  \bleu{ \begin{proposition}
\begin{eqnarray}
    \nabla_{1}(p_n(\mbf{w})) &=& (-1)^{n-1}\,[n]_t[n]_q\, e_n(\mbf{w}),\label{un}\\
    \nabla_{1}(h_n(\mbf{w})) &=&
        \sum_{k=1}^{n} (-q\,t)^{n-k}\, S_{k,1^{n-k}}(\mbf{w}),\label{deux}\\
    \nabla_{1}(e_n(\mbf{w})) &=&
        \sum_{k=1}^{n} [k]_{q,t}e_{n-k}(\mbf{w})e_k(\mbf{w}).\label{trois}
 \end{eqnarray}
  \end{proposition}}
 \begin{proof}[\bf Proof.]\ 
 These are all derived using the fact that
     \begin{equation}
         \nabla_1=\frac{1}{(1-t)(1-q)}\,(\mathrm{Id}-{\mathbf D}_0),
      \end{equation}
   which is directly derived from formula (1.11) of \cite{explicit}. The first one 
   is easiest since we have
\begin{eqnarray*}
    (\mathrm{Id}-{\mathbf D}_0)(p_n(\mbf{w})) &=&
         p_n(\mbf{w})-p_n \!\left[\mbf{w} + \frac{(1-t)(1-q)}{\xi}\right]  
         \Omega'(\mbf{w};-\xi) \big|_{\xi^0}\\
       &=&-(1-t^n)(1-q^n)\frac{1}{\xi^n}  \Omega'(\mbf{w};-\xi)\big|_{\xi^0}\\
       &=&(-1)^{n-1}\,(1-t^n)(1-q^n)\,e_n(\mbf{w}).
\end{eqnarray*}  
For (\ref{deux}), we use
     $$ h_n(\mbf{w}+{\mathbf s})=\sum_{k=0}^n h_k(\mbf{w}) h_{n-k}({\mathbf s})   $$
 to calculate that
 \begin{eqnarray*}
     h_n \!\left[\mbf{w} + \frac{(1-t)(1-q)}{\xi}\right] &=&
         \sum_{k=0}^n \frac{1}{\xi^k} h_{n-k}(\mbf{w} ) h_k((1-t)(1-q))\\
         &=&h_n(\mbf{w})+\sum_{k=1}^n \frac{1}{\xi^k} h_{n-k}(\mbf{w} )  \frac{(1-t)(1-q)(1-q^kt^k)}{(1-q\,t)}.
 \end{eqnarray*}
 Thus we have
     \begin{eqnarray*}
          \nabla_1(h_n(\mbf{w}))&=&\frac{1}{(1-t)(1-q)}\,(\mathrm{Id}-{\mathbf D}_0)(h_n(\mbf{w}))\\
                                &=&-\sum_{k=1}^n h_{n-k}(\mbf{w} )  \frac{(1-q^kt^k)}{(1-q\,t)}\,
                                           \frac{1}{\xi^k}   \Omega'(\mbf{w};-\xi)\big|_{\xi^0}\\
                                &=&-\sum_{k=1}^n (-1)^k  \frac{(1-q^kt^k)}{(1-q\,t)}\,h_{n-k}(\mbf{w})\,e_k(\mbf{w})\\
                                &=&   \sum_{k=1}^{n} (-q\,t)^{n-k}\, S_{k,1^{n-k}}(\mbf{w}).                                    
       \end{eqnarray*}
To get (\ref{trois}), we calculate (now using $e_n(\mbf{w}+{\mathbf s})=\sum_{k=0}^n e_k(\mbf{w}) e_{n-k}({\mathbf s})$) that 
\begin{eqnarray*}
    (\mathrm{Id}-{\mathbf D}_0)(e_n(\mbf{w})) &=&
         e_n(\mbf{w})-e_n \!\left[\mbf{w} + \frac{(1-t)(1-q)}{\xi}\right]  
         \Omega'(\mbf{w};-\xi) \big|_{\xi^0}\\
       &=&-\sum_{k=1}^n \frac{1}{\xi^k} e_k((1-t)(1-q)) e_{n-k}(\mbf{w}) \Omega'(\mbf{w};-\xi)\big|_{\xi^0}\\
       &=&(1-t)(1-q)\sum_{k=1}^n\,[k]_{q,t}\,e_{k}(\mbf{w})e_{n-k}(\mbf{w}).
\end{eqnarray*}
\vskip-30pt
\end{proof}  
We can also see that
\begin{eqnarray}
      \bleu{  \nabla_{n-1}(e_n(\mbf{w})) }&=&   \bleu{\frac{(-1)^{n-1}}{q^{n-1}t^{n-1}}\nabla \nabla_1(h_n(\mbf{w}))}\label{aussi_haglund}\\
                                  &=&   \bleu{ (-1)^{k-1}}
                                          \bleu{    \sum_{k=1}^{n} (-q\,t)^{1-k}\, \nabla(S_{k,1^{n-k}}(\mbf{w}))},
 \end{eqnarray} 
with the second equality obtained using (\ref{deux}).
 Observe that, because of (\ref{epsilon1}), this is also equal to $\nabla_1 \nabla(\varepsilon_{n,1}(\mbf{w};q)/t^{n-1})$ so that we have a link here with (\ref{haglund_epsilon}).
  
  For the symmetric functions $\varepsilon_{n,j}(\mbf{w};q)$, it has been shown (see \cite[(3.34)]{haglund}) that 
         \begin{displaymath}
               q^{\binom{n}{2}}\, \varepsilon_{n,n}(\mbf{w};q)=H_n(\mbf{w};q,t).
          \end{displaymath}
   As already discussed (see (\ref{cocharge_hook})),  it is well-known that $H_n(\mbf{w};q,t)$ expands as
       $$H_n(\mbf{w};q,t)=\sum_{\mu\vdash n} K_{\mu,(n)}(q)\, S_\mu(\mbf{w}),$$
  with the $K_{\mu,(n)}(q)$ polynomials in only the parameter $q$ and having nonnegative integer coefficients.
  More generally, Garsia and Haglund have conjectured that  (see \cite[Conj. 3.8]{haglund}) for all $j$.
 \begin{conj_autres}[Garsia-Haiman 2001]
  \begin{equation}\label{nabla_epsilon}
      \nabla(\varepsilon_{n,j}(\mbf{w};q)) =\sum_{\mu\vdash n} c_\mu^{(n,j)}(q,t)\ S_\mu(\mbf{w}),
  \end{equation}
 with the coefficients $c_\mu^{(n,j)}(q,t)$ polynomials in $q$ and $t$ with positive integer coefficients.
 \end{conj_autres} 
 Moreover, a combinatorial formula is proposed for (\ref{nabla_epsilon}) (see \cite[(6.61)]{haglund}).
  In view of Formula~(\ref{epsilon1}), both (\ref{nabla_epsilon}) and conjecture~\ref{nabla_conj} can be applied here. Sure enough they agree. Formula (\ref{formule_nabla})  implies that we can give a representation theoretic meaning to $\nabla(\varepsilon_{n,n-1}(\mbf{w};q))$. Indeed, we first show that
    \begin{equation}
       \bleu{  \nabla(\varepsilon_{n,n-1}(\mbf{w};q))= {\rm coefficient\ of\ }t\ {\rm in}\ \nabla(e_n(\mbf{w}))},
    \end{equation}
and then use Formula~(\ref{epsilon_n_moins_1}) to calculate $\nabla(\varepsilon_{n,n-1}(\mbf{w};q))$ using results of \cite{science_fiction}.    
Further interest in this arises from the known combinatorial formula \cite[Prop. 3.9.1]{haglund} for the left hand side of this equality.


\end{document}